\documentclass{amsart}
\usepackage{amssymb}

\newtheorem{theorem}{Theorem}[section]

\newtheorem{proposition}[theorem]{Proposition}

\theoremstyle{definition}
\newtheorem{definition}[theorem]{Definition}

\theoremstyle{remark}
\newtheorem{remark}[theorem]{Remark}

\begin{document}

\title{Non-commutative metric topology on matrix state space}
\author{Wei Wu}
\thanks{This research was partially supported by China Scholarship Council, National Natural Science 
Foundation of China, and Shanghai Priority Academic Discipline.}
\address{Department of Mathematics, East China Normal
University, Shanghai 200062, P.R. China}
\email{wwu@math.ecnu.edu.cn} 

\maketitle

\begin{abstract}
We present an operator space version of Rieffel's theorem on the
agreement of the metric topology, on a subset of the Banach space
dual of a normed space, from a seminorm with the weak*-topology.
As an application we obtain a necessary and sufficient condition
for the matrix metric from an unbounded Fredholm module to give
the BW-topology on the matrix state space of the $C^*$-algebra.
Motivated by recent results we formulate a non-commutative
Lipschitz seminorm on a matrix order unit space and characterize
those matrix Lipschitz seminorms whose matrix metric topology
coincides with the BW-topology on the matrix state space.
\end{abstract}

\section{Introduction}\label{sec:1}

One basic idea in modern analysis is that $C^*$-algebras are
non-commutative $C(K)$ spaces\cite{ef,ka}. The basic idea of
non-commutative functional analysis is to study operator spaces as
a generalization of Banach spaces\cite{e,ru,efru}. Let $\mathcal
A$ be a unital $C^*$-algebra. Connes has shown that an appropriate
way to specify a Riemannian metric in this non-commutative
situation is by means of a spectral triple\cite{co}. This consists
of a represention of $\mathcal A$ on a Hilbert space $\mathcal H$,
together with the generalized Dirac operator $D$ on $\mathcal H$
satisfying certain conditions. The Lipschitz seminorm, $L$, is
defined on the set $\mathcal{L(A)}$ of Lipschitz elements of
$\mathcal A$ by the operator norm $L(a)=\|[D, a]\|$. From this,
Connes defines a metric $\rho$ on the state space of $\mathcal A$
by
\[\rho(\mu, \nu)=\sup\{|\mu(a)-\nu(a)|:\,a\in\mathcal{L(A)},
L(a)\le 1\}.\] A natural question is when the metric topology from
the metric $\rho$ agrees with the underlying weak*-topology on the
state space of $\mathcal A$.

In \cite{ri}, Rieffel defined a metric on a subset of the Banach
space dual in a very rudimentary Banach space setting and gave
necessary and sufficient conditions under which the topology from
the metric agrees with the weak*-topology. This problem has also
been considered by Pavlovi\'{c}\cite{pav}. For an operator with
domain contained in a unital $C^*$-algebra, range contained in a
normed space and kernel $\mathbb{C}1$, She defined a metric on the
state space of the $C^*$-algebra from the operator, and then
obtained the necessary and sufficient conditions for the metric to
give the weak*-topology. The theory of operator spaces naturally
leads to the problem of whether we may replace the Banach space by
an operator space. The main purpose of this paper is to extend
these results to the category of operator spaces.

The paper is organized as follows. We begin in section \ref{sec:2}
with the explanation of terminology and notation. In section
\ref{sec:3} we present the operator space version of Rieffel's
theorem and in section \ref{sec:4} we use it to the matrix metrics
from Dirac operators. Finally in section \ref{sec:5} we define a
non-commutative version of Lipschitz seminorms on matrix order
unit spaces and then use the results of the previous sections to
determine when the metric topology from a matrix Lipschitz
seminorm agrees with the BW-topology on the matrix state space.

\section{Terminology and notation}\label{sec:2}

All vector spaces are assumed to be complex throughout this paper.
Given a vector space $V$, we let $M_{m, n}(V)$ denote the matrix
space of all $m$ by $n$ matrices $v=[v_{ij}]$ with $v_{ij}\in V$,
and we set $M_n(V)=M_{n,n}(V)$. If $V=\mathbb C$,  we write
$M_{m,n}=M_{m,n}(\mathbb C)$ and $M_n=M_{n,n}(\mathbb C)$, which
means that we may identify $M_{m,n}(V)$ with the tensor product
$M_{m,n}\otimes V$. We identify $M_{m,n}$ with the normed space
$\mathcal B({\mathbb C}^n, {\mathbb C}^m)$. We use the standard
matrix multiplication and *-operation for compatible scalar
matrices, and $1_n$ for the identity matrix in $M_n$. There are
two natural operations on the matrix spaces. For $v\in M_{m,n}(V)$
and $w\in M_{p,q}(V)$, the direct sum $v\oplus w\in
M_{m+p,n+q}(V)$ is defined by letting $v\oplus
w=\left[\begin{array}{cc}v&0\\ 0&w\end{array}\right]$, and if we
are given $\alpha\in M_{m,p}$, $v\in M_{p,q}(V)$ and $\beta\in
M_{q,n}$, the matrix product  $\alpha v\beta\in M_{m,n}(V)$ is
defined by $\alpha
v\beta=\left[\sum_{k,l}\alpha_{ik}v_{kl}\beta_{lj}\right]$.

If $V$ and $W$ are vector spaces in duality, then it determines
the matrix pairing $\ll\cdot, \cdot\gg:\, M_n(V)\times
M_m(W)\longmapsto M_{nm}$, where $\ll[v_{ij}],
[w_{kl}]\gg=\left[<v_{ij}, w_{kl}>\right]$ for $[v_{ij}]\in
M_n(V)$ and $[w_{kl}]\in M_m(W)$. Given an arbitrary vector space
$V$, a {\it matrix gauge} $\mathcal G=(g_n)$ on $V$ is a sequence
of gauges $g_n:\, M_n(V)\longmapsto [0, +\infty]$ such that
$g_{m+n}(v\oplus w)=\max\{g_m(v), g_n(w)\}$ and $g_n(\alpha
v\beta)\le\|\alpha\| g_m(v)\|\beta\|$ for any $v\in M_m(V)$, $w\in
M_n(V)$, $\alpha\in M_{n,m}$ and $\beta\in M_{m,n}$. A matrix
gauge $\mathcal G=(g_n)$ is a {\it matrix seminorm} on $V$ if for
any $n\in\mathbb N, g_n(v)<+\infty$ for all $v\in M_n(V)$. If each
$g_n$ is a norm on $M_n(V)$, we say that $\mathcal{G}$ is a {\it
matrix norm}. An {\it operator space} is a vector space together
with a matrix norm. For other terminology and notation, we will
follow \cite{chef} and \cite{efwe}.

\section{Matrix metrics on matrix states}\label{sec:3}

Let $V$ be an operator space with the matrix norm
$\|\cdot\|=(\|\cdot\|_n)$, and let $W$ be a subspace of $V$. As in \cite{ri}, 
we will not assume that $W$ is closed. Suppose that $\mathcal{L}=(L_n)$ is a matrix seminorm on $W$ and 
$\varphi_1^{(0)}$ is a bounded linear functional on the subspace 
$K_1=\{a\in W:\, L_1(a)=0\}$ with $\|\varphi^{(0)}_1\|=1$. 
Set $\varphi^{(0)}_n=\underbrace{\varphi^{(0)}_1\oplus\cdots\oplus\varphi^{(0)}_1}_n$, $n\in\mathbb{N}$. Then
$(\varphi^{(0)}_n)$ is a sequence of completely bounded linear
mappings $\varphi^{(0)}_n:\,K_1\mapsto M_n$, with $\|\varphi^{(0)}_n\|_{cb}=1$. Here 
$\|\varphi^{(0)}_n\|_{cb}$ is the completely bounded norm of
$\varphi^{(0)}_n$. In the applications, $K_1$ will always be
closed in $V$. Thus we assume that $K_1$ is closed in $V$. Denote
$K_n=\{a\in M_n(W): L_n(a)=0\}$ for $n\in\mathbb N$, and let
$\mathcal{K}=(K_n)$. The inequality $L_1(v_{ij})\le
L_n([v_{ij}])\le\sum_{i,j=1}^nL_1(v_{ij})$ holds for $[v_{ij}]\in
M_n(W)$ and $n\in\mathbb N$\cite{efwe}. So $K_n=M_n(K_1)$ for
$n\in\mathbb N$.

Denote $GCS_n(V)=\{\varphi:\, \varphi$ is a completely bounded
linear mapping from $V$ to $M_n$, $\varphi=\varphi_n^{(0)}$ on
$K_1$, and $\|\varphi\|_{cb}=1\}$, and set $\mathcal
{GCS}(V)=(GCS_n(V))$. We write $CB_n(V)=\{\varphi:\, \varphi$ is a
completely bounded linear mapping from $V$ to $M_n\}$, and
$CB_n^1(V)=\{\varphi\in CB_n(V):\, \|\varphi\|_{cb}\le 1\}$, for $n\in\mathbb N$. 
Then $GCS_n(V)$ is a convex subset of $CB_n(V)$. Define
$S:\,M_n(V^*)\mapsto CB_n(V)$ by
$(S([\rho_{ij}]))(a)=[\rho_{ij}(a)]$, where $[\rho_{ij}]\in
M_n(V^*)$ and $a\in V$. Then $S$ is a linear isomorphism from
$M_n(V^*)$ onto $CB_n(V)$. As a linear mapping from operator space
$M_n(V^*)$ onto operator space $CB_n(V)$, $S$ is a completely
isometry\cite{bl, efr}. On $\mathcal{GCS}(V)$, we can define a
sequence $\mathcal{D_L}=(D_{L_n})$ of metrics by
\[D_{L_n}(\varphi, \psi)=\sup\{\|\ll\varphi, a\gg-\ll\psi,
a\gg\|:\, a\in M_r(W), L_r(a)\le 1, r\in\mathbb N\},\]for
$\varphi, \psi\in GCS_n(V)$. (These may take value $+\infty$)We
call $\mathcal{D_L}$ the {\it matrix metric} from $\mathcal{L}$.
We will refer to the topology on $\mathcal {GCS}(V)$ defined by
$\mathcal{D_L}$, that is, the topology on each $GCS_n(V)$ is
induced by $D_{L_n}$, as the $\mathcal{D_L}$-{\it topology} or the
{\it matrix metric topology}. The natural topology on
$\mathcal{GCS}(V)$ is the {\it BW-topology}, that is, topologies
each $GCS_n(V)$ by BW-topology\cite{ar, pau}. We say that $W$ {\it
separates the points} of $\mathcal{GCS}(V)$ if given $\varphi,
\psi\in GCS_n(V)$ and $\varphi\neq\psi$ for $n\in\mathbb{N}$ there
is an $a\in W$ such that $\varphi(a)\neq\psi(a)$.

\begin{proposition}\label{pro:301}
For $\varphi, \psi\in GCS_n(V)$,
\[D_{L_n}(\varphi, \psi)=\sup\{\|\ll \varphi, a\gg-\ll\psi,
a\gg\|: a\in M_n(W), L_n(a)\le 1\}.\]
\end{proposition}

\begin{proof}
Given $\varphi, \psi\in GCS_n(V)$. Let $K_{L_n}(\varphi,
\psi)=\sup\{\|\ll \varphi, a\gg-\ll\psi, a\gg\|: a\in M_n(W),
L_n(a)\le 1\}$. Clearly, $K_{L_n}(\varphi, \psi)\le
D_{L_n}(\varphi, \psi)$. For any $a\in M_n(W)$, we have that
$\|\ll\varphi, a\gg-\ll\psi, a\gg\|\le K_{L_n}(\varphi,
\psi)L_n(a)$. Since $\mathcal{L}$ is a matrix gauge on $W$, we
have that $\|\ll\varphi, a\gg-\ll\psi, a\gg\|\le K_{L_n}(\varphi,
\psi)L_r(a)$ for all $a\in M_r(W)$ with $r\in\mathbb{N}$ arbitrary
by Lemma 5.2 in \cite{efwe}. So $D_{L_n}(\varphi,
\psi)=\sup\{\|\ll \varphi, a\gg-\ll\psi, a\gg\|: a\in M_r(W),
L_r(a)\le 1, r\in\mathbb{N}\}\le K_{L_n}(\varphi, \psi)$.
Therefore, $K_{L_n}=D_{L_n}$.
\end{proof}

\begin{proposition}\label{pro:31} $\mathcal{GCS}(V)$ is BW-compact.
\end{proposition}

\begin{proof} This follows immediately from the observation that $GCS_n(V)=
\bigcap_{k\in K_1}\{\varphi\in CB_n^1(V):\,\varphi(k)=
\varphi_n^{(0)}(k)\}$ and the fact that $CB_n^1(V)$ is
BW-compact\cite{pau}.
\end{proof}

\begin{proposition}\label{pro:32} If $W$ separates the points of $\mathcal{GCS}(V)$,
then the $\mathcal{D_L}$-topology on $\mathcal {GCS}(V)$ is finer
than the BW-topology.
\end{proposition}

\begin{proof}
Let $\{\varphi_k\}$ be a sequence in $GCS_n(V)$ which converges to
$\varphi\in GCS_n(V)$ for the metric $D_{L_n}$. Then for any
$r\in\mathbb N$ and $a\in M_r(W)$ with $L_r(a)\le 1$, we have
$\lim_{k\to\infty}\|\ll \varphi_k, a\gg-\ll\varphi, a\gg\|=0$. So,
for every $r\in\mathbb{N}$ and $a\in M_r(W)$, $\{\ll\varphi_k,
a\gg\}$ converges to $\ll\varphi, a\gg$ in the norm topology, and
hence in the weak operator topology.

Denote $\phi(a, \xi, \eta)(\psi)=<\ll\psi, a\gg\xi, \eta>$ for
$a\in M_r(W)$, $\psi\in GCS_n(V)$, $\xi, \eta\in\mathbb{C}^{rn}$
and $r, n\in\mathbb{N}$. The paragraph above shows that $\phi(a,
\xi, \eta)(\varphi_k)$ converges to $\phi(a, \xi, \eta)(\varphi)$
for any $a\in M_r(W)$, $\xi, \eta\in\mathbb{C}^{rn}$  and
$r\in\mathbb{N}$. Let $Q_{ r,n}=\{\phi(a, \xi, \eta):\,a\in
M_r(W), \xi, \eta\in\mathbb{C}^{rn}\}$ and
$Q_n=\bigcup_{r=1}^\infty Q_{r,n}$. Then $Q_n$ is a linear space
of BW-continuous functions on $GCS_n(V)$. Since $W$ separates the
points of $\mathcal{GCS}(V)$, $Q_n$ separates the points of
$GCS_n(V)$. For $a\in K_r$ and $\xi, \eta\in\mathbb{C}^{rn}$, we
have that $\phi(a, \xi, \eta)(\varphi)=<\ll\varphi, a\gg\xi,
\eta>=<\ll\varphi_n^{(0)}, a\gg\xi, \eta>$, where $\varphi\in
GCS_n(V)$. So $Q_n$ contains the constant functions. $GCS_n(V)$ is
BW-compact implies that $Q_n$ determines the BW-topology of
$GCS_n(V)$. Therefore, $\{\varphi_k\}$ converges to $\varphi$ in
the BW-topology.
\end{proof}

For $a\in M_n(W)$, define $\|a\|_n^\natural=\sup\{\|\ll\varphi,
a\gg\|:\, \varphi\in GCS_n(V)\}$. Then we obtain a sequence
$\|\cdot\|^\natural=(\|\cdot\|_n^\natural)$ of seminorms, and
$\|a\|_n^\natural\le\|a\|_n$ for $a\in M_n(W)$.

\begin{proposition}\label{pro:33} If $\|\cdot\|^\natural$ is a matrix
seminorm on $W$, then
\[\|a\|_n^\natural=\sup\{\|\ll\varphi,
a\gg\|:\, \varphi\in GCS_r(V), r\in\mathbb N\}.\]
\end{proposition}

\begin{proof} Denote
$\|a\|_n^\flat=\sup\{\|\ll\varphi, a\gg\|:\, \varphi\in GCS_r(V),
r\in\mathbb N\}$. Clearly, $\|a\|_n^\natural\le\|a\|_n^\flat$.
Given $\varphi\in GCS_r(V)$. If $r<n$, we choose a $\varphi_1\in
GCS_{n-r}(V)$ and then obtain that $\|\ll\varphi,
a\gg\|\le\max\{\|\ll\varphi,a\gg\|, \|\ll\varphi_1,a\gg\|\}
=\left\|\left[\begin{array}{cc}
\ll\varphi, a\gg&0\\
0&\ll\varphi_1,a\gg\end{array}\right]\right\|
=\|\ll\varphi\oplus\varphi_1, a\gg\|\le\|a\|_n^\natural$. When
$r>n$, we have that $\|\ll\varphi, a\gg\|=\|\ll\varphi, a\oplus
0_{r-n}\gg\|\le\|a\oplus 0_{r-n}\|_r^\natural=\|a\|_n^\natural$
since $\|\cdot\|^\natural$ is a matrix seminorm. Hence
$\|a\|_n^\flat\le\|a\|_n^\natural$, and this, with the previous
inequality, gives that $\|a\|_n^\natural=\sup\{\|\ll\varphi,
a\gg\|:\, \varphi\in GCS_r(V), r\in\mathbb N\}$.
\end{proof}

Let $\tilde{W}=W/K_1$ and $\tilde{V}=V/K_1$. Then $\mathcal{L}$
drops to a matrix norm on $\tilde{W}$\cite{efwe}, which we denote
by $\tilde{\mathcal{L}}=(\tilde{L}_n)$. But on $\tilde{W}$ and
$\tilde{V}$ we also have the quotient matrix norm from $\|\cdot\|$
on $V$\cite{efru}, which we denote by
$\|\cdot\|^\sim=(\|\cdot\|_n^\sim)$. The image in $M_n(\tilde{W})$
of $a\in M_n(W)$ will be denoted by $\tilde{a}$. For each
$n\in\mathbb N$, we let $L_n^t=\{a\in M_n(W):\, L_n(a)\le t\}$ and
$B_n^t=\{a\in M_n(W): L_n(a)\le 1, \|a\|_n\le t\}$ for $t>0$, and
set $\mathcal{L}^t=(L_n^t)$ and $\mathcal{B}^{t}=(B_n^t)$. We say
that the image of a graded set $\mathbf{E}=(E_n)$($E_n\subseteq
M_n(W)$) in $\tilde{\mathcal W}$ is {\it totally bounded} for
$\|\cdot\|^\sim$ if each $\tilde{E}_n$ is totally bounded in
$M_n(\tilde{\mathcal W})$ for $\|\cdot\|^\sim_n$.

\begin{definition}\label{def:34} Let $V$ be a vector space, and let $\|\cdot\|=(\|\cdot\|_n)$ and
$\|\cdot\|^\sharp=(\|\cdot\|_n^\sharp)$ be two matrix norms on
$V$. We say that $\|\cdot\|$ and $\|\cdot\|^\sharp$ are {\it
equivalent} if for any $n\in\mathbb{N}$, $\|\cdot\|_n$ is
equivalent to $\|\cdot\|_n^\sharp$.
\end{definition}

\begin{proposition}\label{pro:35}
Suppose that $\|\cdot\|=(\|\cdot\|_n)$ and
$\|\cdot\|^\sharp=(\|\cdot\|_n^\sharp)$ be two matrix norms on the
vector space $V$. If $\|\cdot\|_1$ and $\|\cdot\|_1^\sharp$ are
equivalent, then $\|\cdot\|$ and $\|\cdot\|^\sharp$ are
equivalent.
\end{proposition}

\begin{proof}
Since $\|\cdot\|_1$ and $\|\cdot\|_1^\sharp$ are equivalent, there
are positive constants $\alpha$ and $\beta$ such that
$\alpha\|x\|_1^\sharp\le\|x\|_1\le\beta\|x\|_1^\sharp$ for $x\in
V$. Because $\|\cdot\|$ and $\|\cdot\|^\sharp$ are matrix norms,
we have
$\|x_{ij}\|_1\le\|[x_{ij}]\|_n\le\sum_{i,j=1}^n\|x_{ij}\|_1$ and
$\|x_{ij}\|_1^\sharp\le\|[x_{ij}]\|_n^\sharp\le\sum_{i,j=1}^n\|x_{ij}
\|_1^\sharp$ for $[x_{ij}]\in M_n(V)$\cite{ru}. Hence $\alpha
n^{-2}\|[x_{ij}]\|_n^\sharp\le\|[x_{ij}]\|_n\le
n^2\beta\|[x_{ij}]\|_n^\sharp$. So $\|\cdot\|_n$ and
$\|\cdot\|_n^\sharp$ are equivalent. Therefore, $\|\cdot\|$ and
$\|\cdot\|^\sharp$ are equivalent.
\end{proof}

\begin{theorem}\label{th:36} Let $V$ be an operator space with the matrix norm
$\|\cdot\|=(\|\cdot\|_n)$, and let $W$ be a subspace of $V$, which is not necessarily closed. 
Suppose that $\mathcal{L}=(L_n)$ is a matrix seminorm on $W$ with
closed $($in $V$$)$ kernels $K_n=\{a\in M_n(W):\,
L_n(a)=0\}(n\in\mathbb{N})$ and $(\varphi^{(0)}_n)$ is a sequence
of completely bounded linear mappings
$\varphi^{(0)}_n:\,K_1\mapsto M_n$, with 
$\varphi^{(0)}_n=\underbrace{\varphi^{(0)}_1\oplus\cdots\oplus\varphi^{(0)}_1}_n$ and 
$\|\varphi^{(0)}_n\|_{cb}=1$. Denote $GCS_n(V)=\{\varphi: \varphi$
is a completely bounded linear mapping from $V$ to $M_n$,
$\varphi=\varphi_n^{(0)}$ on $K_1$, and $\|\varphi\|_{cb}=1\}$ and
$\|a\|_n^\natural=\sup\{\|\ll\varphi, a\gg\|: \varphi\in
GCS_n(V)\}, a\in M_n(W)$. Set $\mathcal {GCS}(V)=(GCS_n(V))$ and
$\|\cdot\|^\natural=(\|\cdot\|^\natural_n)$. We assume that $W$
separates the points of $\mathcal{GCS}(V)$.
\begin{enumerate}
\item Suppose that $\|\cdot\|^\natural$ is a matrix norm on $W$ and it is
equivalent to the matrix norm $\|\cdot\|$. If the
$\mathcal{D_L}$-topology on $\mathcal{GCS}(V)$ agrees with the
BW-topology, then the image of $\mathcal{L}^1$ in $\tilde{W}$ is
totally bounded for $\|\cdot\|^\sim$.

\item If the image of $L_1^1$ in $\tilde{W}$ is totally bounded
for $\|\cdot\|_1^\sim$, then the $\mathcal{D_L}$-topology on
$\mathcal{GCS}(V)$ agrees with the BW-topology.
\end{enumerate}
\end{theorem}

\begin{proof} (1) If $\|\cdot\|^\natural$ is
equivalent to $\|\cdot\|$, then $\|\cdot\|_1^\natural$ is
equivalent to $\|\cdot\|_1$. When $n=1$, $GCS_1(V)$ is just the
space $S$ in \cite{ri}. Since the $D_{L_1}$-topology agrees with
the BW-topology on $GCS_1(V)$, the image of $L_1^1$ in $\tilde{W}$
is totally bounded for $\|\cdot\|_1^\sim$ by Theorem 1.8 in
\cite{ri} and Proposition \ref{pro:301}. So given $\epsilon>0$,
there exist elements $a_1, \cdots,a_k\in L_1^1$ such that
$\tilde{L}_1^1\subseteq\bigcup_{j=1}^k\left\{\tilde{a}\in
\tilde{W}: \|\tilde{a}-\tilde{a}_j\|_1^\sim<\epsilon\right\}$. For
a matrix gauge $\mathcal{G}=(g_n)$ on a vector space $Y$ and
$v=[v_{ij}]\in M_n(Y)$, we have the constraint $g_1(v_{ij})\le
g_n(v)\le\sum_{i,j=1}^ng_1(v_{ij})$ \cite{efwe}. From the first
inequality, we get that $L_n^1\subseteq M_n(L_1^1)$.

Denote $X_n=\{a\in M_n(W): a=[a_{ij}],
a_{ij}\in\{a_1,\cdots,a_k\}\}$. Then $X_n$ is a finite set for
every $n\in\mathbb N$. Now for $b=[b_{ij}]\in L_n^1$, there exists
an $a=[a_{ij}]\in X_n$ such that
$\|\tilde{b}_{ij}-\tilde{a}_{ij}\|_1^\sim<\epsilon$. So
$\|\tilde{b}-\tilde{a}\|_n^\sim\le\sum_{i,j=1}^n\|\tilde{b}_{ij}
-\tilde{a}_{ij}\|_1^\sim<n^2\epsilon$. Therefore, the image of
$L_n^1$ in $M_n(\tilde{W})$ is totally bounded for
$\|\cdot\|_n^\sim$. Consequently, the image of $\mathcal{L}^1$ in
$\tilde{W}$ is totally bounded for $\|\cdot\|^\sim$.

(2) Suppose that the image of $L_1^1$ in $\tilde{W}$ is totally
bounded for $\|\cdot\|_1^\sim$. From the proof of (1), we see that
the image of $\mathcal{L}^1$ in $\tilde{W}$ is totally bounded for
$\|\cdot\|^\sim$. For any $n\in\mathbb N$, let $\varphi\in
GCS_n(V)$ and $\epsilon>0$ be given, and let $B(\varphi,
\epsilon)$ be the $D_{L_n}$-ball of radius $\epsilon$ at
$\varphi$. By the total boundedness of $\tilde{L}_n^1$, we can
find elements $a_1, a_2, \cdots, a_k\in L_n^1$ such that the
$\|\cdot\|^\sim_n$-balls of radius $\frac{\epsilon}6$ at
$\tilde{a}_j$'s cover $\tilde{L}^1_n$.

Assume that $\varphi=[\varphi_{pq}]$ and $a_j=[a_{st}^{(j)}]$ for
$j\in\{1,2,\cdots,k\}$. Let $Q=\big\{\psi=[\psi_{ij}]\in GCS_n(V):
|\psi_{pq}(a_{st}^{(j)})-\varphi_{pq}(a_{st}^{(j)})|<
\frac{\epsilon}{6n^4},
s,t,p,q\in\{1,2,\cdots,n\},j\in\{1,2,\cdots,k\}\big\}$. Then $Q$
is BW-open and $\varphi\in Q$. Set $N=\{\psi\in GCS_n(V):\,
\|\ll\varphi-\psi, a_j\gg\|<\frac{\epsilon}6, 1\le j\le k\}$. Then
$Q\subseteq N$ and so $N$ is a BW-neighborhood of $\varphi$. For
any $a\in L_n^1$, there are a $j\in\{1,2,\cdots,k\}$ and a $b\in
K_n$ such that $\|a-(a_j+b)\|_n<\frac{\epsilon}6$. Hence for any
$\psi\in N$, we have that $\|\ll\varphi-\psi,
a\gg\|\le\|\varphi\|_{cb}\|a-(a_j+b)\|_n+\|\ll\varphi-\psi,
a_j\gg\|+ \|\psi\|_{cb}\|(a_j+b)-a\|_n<\frac{\epsilon}{2}$. By
Lemma 5.2 in \cite{efwe}, we obtain $\|\ll\varphi-\psi,
a\gg\|\le\frac{\epsilon}{2}$ for all $a\in L_r^1$ with
$r\in\mathbb N$ arbitrary. Thus $D_{L_n}(\varphi,
\psi)\le\frac{\epsilon}{2}$. Consequently $N\subseteq B(\varphi,
\epsilon)$. Since $W$ separates the points of $\mathcal{GCS}(V)$,
by Proposition \ref{pro:32} the $D_{L_n}$-topology on $GCS_n(V)$
is finer that the BW-topology. Therefore, the $D_{L_n}$-topology
on $GCS_n(V)$ agrees with the BW-topology. We conclude that the
$\mathcal{D_L}$-topology on $\mathcal{GCS}(V)$ agrees with the
BW-topology.
\end{proof}

For a sequence $\{r_n\}$ of positive constants, we set
$CB_n^{r_n}(\tilde{V})=\{f\in M_n((\tilde{V})^*): \|f\|_{cb}\le
r_n\}$, and this is just the subset of $M_n(V^*)$ consisting of
those $f\in M_n(V^*)$ such that $\|f\|_{cb}\le
r_n$ and $f(k)=0_n$ for $k\in K_1$. If
$CB_n^{r_n}(\tilde{V})\subseteq
(GCS_n(V)-GCS_n(V))+i(GCS_n(V)-GCS_n(V))$ for each
$n\in\mathbb{N}$, we say that $V$ has the {\it matrix
$\{r_n\}$-decomposable property} about $\mathcal{K}$; if $r_n=r$
for all $n\in\mathbb{N}$, we say that $V$ has the {\it matrix
$r$-decomposable property} about $\mathcal{K}$.

\begin{proposition}\label{pro:37}
If $V$ has the matrix $\{r_n\}$-decomposable property about
$\mathcal{K}$ and $\mathcal{D_L}$ is bounded, that is, each
$D_{L_n}$ is bounded, then there is a sequence $\{d_n\}$ of
positive constants such that $\|\tilde{a}\|^\sim_n\le d_n$ for all
$a\in L_n^1$ and $n\in\mathbb{N}$.
\end{proposition}

\begin{proof}
Assume that each $D_{L_n}$ is bounded by $c_n$. Then for $a\in
M_n(W)$ with $L_n(a)\le 1$, we have that
$\|\tilde{a}\|_n^\sim=\sup\{\|\ll f, \tilde{a}\gg\|: f\in
CB_n(\tilde{V}), \|f\|_{cb}\le 1\}\le r_n^{-1}\sup\{\|\ll
(\varphi_1-\varphi_2)+i(\varphi_3-\varphi_4), \tilde{a}\gg\|:
\varphi_i\in GCS_n(V), i=1, 2, 3, 4\}\le 2r_n^{-1}\sup\{\|\ll
\varphi_1-\varphi_2, \tilde{a}\gg\|: \varphi_1, \varphi_2\in
GCS_n(V)\}=2r_n^{-1}\sup\{\|\ll \varphi_1, a\gg-\ll\varphi_2,
a\gg\|: \varphi_1, \varphi_2\in GCS_n(V)\} \le
2r_n^{-1}\sup\{D_{L_n}(\varphi_1, \varphi_2): \varphi_1,
\varphi_2\in GCS_n(V)\}\le 2c_nr_n^{-1}$ by Theorem 4.2 in
\cite{ru} and Lemma 2.1 in \cite{efr}. Letting $d_n=2c_nr_n^{-1}$,
we have that $\|\tilde{a}\|_n^\sim\le d_n$ for all $a\in L_n^1$.
\end{proof}

Here, we give another additional premise such that for the
$\mathcal{D_L}$-topology on $\mathcal{GCS}(V)$ to agree with the
BW-topology it is necessary that the image of $\mathcal{L}^1$ in
$\tilde{W}$ is totally bounded for $\|\cdot\|^\sim$.

\begin{theorem}\label{th:38}
With the notation of Theorem \ref{th:36}, suppose that $V$ has the
matrix $\{r_n\}$-decomposable property about $\mathcal{K}=(K_n)$.
If the $\mathcal{D_L}$-topology on $\mathcal{GCS}(V)$ agrees with
the BW-topology, then the image of $\mathcal{L}^1$ in $\tilde{W}$
is totally bounded for $\|\cdot\|^\sim$.
\end{theorem}

\begin{proof}
If each $D_{L_n}$-topology agrees with the BW-topology on
$GCS_n(V)$, then each $D_{L_n}$ must be bounded since $GCS_n(V)$
is BW-compact. By Proposition \ref{pro:37}, there exists a
sequence $\{d_n\}$ of positive constants such that
$\|\tilde{a}\|_n^\sim\le d_n$ for $a\in L_n^1$ and
$n\in\mathbb{N}$. Choose $k_n>d_n$. Then
$\|\tilde{a}\|_n^\sim<k_n$ if $a\in L_n^1$. So the image of
$B_n^{k_n}$ in $M_n(\tilde{W})$ is the same as the image of
$L_n^1$. For $r,n\in\mathbb N$, $a\in L_n^1$ and $\varphi, \psi\in
GCS_r(V)$, we have that $\|\ll\varphi, a\gg-\ll\psi, a\gg\|\le
D_{L_r}(\varphi, \psi)$. Let $\epsilon>0$. For each $\varphi\in
GCS_n(V)$, let $U_\varphi$ be an open neighborhood of $\varphi$
such that $\|\ll\varphi, a\gg-\ll\psi, a\gg\|<\frac{\epsilon}{3}$
for $a\in L_n^1(n\in\mathbb{N})$, and $\psi\in U_\varphi$. Now
$\{U_\varphi:\,\varphi\in GCS_n(V)\}$ is an open covering of
$GCS_n(V)$. Since $GCS_n(V)$ is BW-compact, there are points
$\varphi_1, \cdots, \varphi_k$ in $GCS_n(V)$ such that
$GCS_n(V)=\cup_{j=1}^kU_{\varphi_j}$. Assume that
$\varphi_j=[\varphi^{(j)}_{st}], j=1, \cdots, k$, and denote $D =
\{ \ddot{a}: \ddot{a}$ $= (\varphi^{(1)}_{11} (a_{11}),
\varphi^{(1)}_{11} (a_{12})$, $\cdots, \varphi^{(1)}_{11}
(a_{nn}), \varphi^{(1)}_{12} (a_{11}), \varphi^{(1)}_{12}
(a_{12}), \cdots, \varphi_{12}^{(1)} (a_{nn}), \cdots$,
$\varphi_{nn}^{(1)} (a_{11}), \varphi_{nn}^{(1)} (a_{12}), \cdots,
\varphi^{(1)}_{nn} (a_{nn})$, $\varphi^{(2)}_{11} (a_{11})$,
$\cdots$, $\varphi^{(2)}_{nn} (a_{nn})$, $\cdots$,
$\varphi^{(k)}_{11} (a_{11}), \cdots$, $\varphi^{(k)}_{nn}$
$(a_{nn})), a= [a_{ij}]\in B_n^{k_n}\}$. For $a\in B_n^{k_n}$, we
have
\[\begin{array}{rcl}
\|\ddot{a}\|&=& \sqrt{\sum_{j=1}^k\sum_{s,t=1}^n\sum_{r,l=1}^n
|\varphi^{(j)}_{st}(a_{rl})|^2}\\
&=&\sqrt{\sum_{j=1}^k\sum_{s,t=1}^n\sum_{r,l=1}^n
|e_{sr}^*\ll\varphi_j, a\gg e_{tl}|^2}\\
&\le&\sqrt{\sum_{j=1}^k\sum_{s,t=1}^n\sum_{r,l=1}^n\|\ll\varphi_j,
a\gg\|^2}\\
&\le&\sqrt{\sum_{j=1}^k\sum_{s,t=1}^n\sum_{r,l=1}^n\|\varphi_j\|_{cb}^2\|a\|_n^2}\le
n^2k_n\sqrt{k},\end{array}\] where $e_{ij}$ is the $nr\times 1$
column matrix $[0, \cdots, 0, 1_{(i-1)n+j}, 0, \cdots, 0]^*$. So
$D$ is totally bounded, and hence there are elements $a_1,
a_2,\cdots, a_m\in B_n^{k_n}$ such that
$D\subseteq\bigcup_{j=1}^m\left\{\ddot{a}:
\|\ddot{a}-\ddot{a}_j\|< \frac{\epsilon}{3n^2}, a\in B_n^{k_n}\right\}$.
Now for $a\in B_n^{k_n}$, there exists a $p\in \{1,2,\cdots,m\}$
such that $\|\ddot{a}-\ddot{a}_p\|<\frac{\epsilon}{3}$. For any
$\varphi\in GCS_n(V)$, there is a $j\in\{1,2,\cdots,k\}$ such that
$\varphi\in U_{\varphi_j}$. So by Cauchy's inequality, we have $ \|\ll a-a_p, \varphi\gg\|
\le\|\ll a-a_p, \varphi_j\gg\|+\|\ll a-a_p, \varphi-\varphi_j\gg\|
\le n^2\|\ddot{a}-\ddot{a}_p\|+\|\ll a,
\varphi-\varphi_j\gg\|+\|\ll a_p,
\varphi-\varphi_j\gg\|<\epsilon$. By Lemma 2.1 in \cite{efr}, we
have that
$\|\tilde{a}-\tilde{a}_p\|_n^\sim=\sup\{\|\ll\tilde{a}-\tilde{a}_p,
f\gg\|: f\in M_n(\tilde{V}), \|f\|_{cb}\le 1\}
=r_n^{-1}\sup\{\|\ll\tilde{a}-\tilde{a}_p,
(\varphi_1-\varphi_2)+i(\varphi_3-\varphi_4)\gg\|: \varphi_i\in
GCS_n(V), i=1,2,3,4\} \le
2r_n^{-1}\sup\{\|\ll\tilde{a}-\tilde{a}_p,
\varphi_1-\varphi_2\gg\|: \varphi_1,\varphi_2\in GCS_n(V)\}
=2r_n^{-1}\sup\{\|\ll a-a_p, \varphi_1-\varphi_2\gg\|:
\varphi_1,\varphi_2\in GCS_n(V)\}\le 4r_n^{-1}\sup\{\|\ll a-a_p,
\varphi\gg\|: \varphi\in GCS_n(V)\}\le 4r_n^{-1}\epsilon$.
Therefore, the image of $B_n^{k_n}$ in $M_n(\tilde{W})$ is totally
bounded for $\|\cdot\|_n^\sim$, then so is $L_n^1$. From the
arbitrariness of $n$, the image of $\mathcal{L}^1$ in $\tilde{W}$
is totally bounded for $\|\cdot\|^\sim$.
\end{proof}

\section{Matrix metrics from Dirac operators}\label{sec:4}

Let $\mathcal{A}$ be a unital $C^*$-algebra. An {\it unbounded
Fredholm module} $(\mathcal{H}, D)$ over $\mathcal{A}$
is\cite{co,pav}: a Hilbert space $\mathcal{H}$ which is a left
$\mathcal{A}$-module, that is, a Hilbert space $\mathcal{H}$ and a
*-representation of $\mathcal{A}$ on $\mathcal{H}$; an unbounded,
self-adjoint operator $D$ (the generalized Dirac operator) on
$\mathcal{H}$ such that $\mathcal{B}=\{a\in\mathcal{A}$: $[D, a]$
is densely defined and extends to a bounded operator on
$\mathcal{H}\}$ is norm dense in $\mathcal{A}$; $(1+D^2)^{-1}$ is
a compact operator. Given a triple $(\mathcal{A}, \mathcal{H},
D)$, where $(\mathcal{H}, D)$ is an unbounded Fredholm module over
a unital $C^*$-algebra $\mathcal{A}$, then $\mathcal{B(H)}$ is a
Banach $\mathcal{A}$-module, and the mapping
$T:\,\mathcal{A}\mapsto\mathcal{B(H)}$, defined by $Ta=[D, a]$, is a
densely defined derivation from $\mathcal{A}$ into
$\mathcal{B(H)}$. For $n\in\mathbb{N}$, we define
$L_n(a)=\|T_n(a)\|, a\in M_n(\mathcal{B})$, and denote
$\mathcal{L}=(L_n)$. For the unital $C^*$-algebra $\mathcal{A}$,
the {\it matrix state space} of $\mathcal{A}$ is the collection
$\mathcal{CS(A)}=(CS_n(\mathcal{A}))$ of matrix states
$CS_n(\mathcal{A})=\{\varphi: \varphi$ is a completely positive
linear mapping from $\mathcal{A}$ to $M_n$, $\varphi(1)=1_n\}$,
where $1$ is the identity of $\mathcal{A}$. This space is
important in the study of $C^*$-algebras\cite{smw}.

\begin{proposition}\label{pro:41} $\mathcal{L}$ is a matrix seminorm on
$\mathcal{B}$, and $M_n(\mathbb{C}1)\subseteq K_n=\{a\in
M_n(\mathcal{B}):\,L_n(a)=0\}$.
\end{proposition}

\begin{proof} It is clear that $L_n$ is a seminorm on $M_n(\mathcal{B})$ for each
$n\in\mathbb{N}$. For $a=[a_{ij}]\in M_m(\mathcal{B})$,
$b=[b_{st}]\in M_n(\mathcal{B})$, $\alpha=[\alpha_{kl}]\in
M_{n,m}$ and $\beta=[\beta_{uv}]\in M_{m,n}$, we have that
$L_{m+n}(a\oplus b)=\|T_{m+n}(a\oplus b)\|=\|T_m(a)\oplus T_n(b)\|
=\|[[D, a_{ij}]]\oplus[[D, b_{st}]]\| =\max\{\|[[D, a_{ij}]]\|,
\|[[D, b_{st}]]\|\}=\max\{\|T_m(a)\|, \|T_n(b)\|\}=\max\{L_m(a), L_n(b)\}$ and 
$L_n(\alpha a\beta)=\|T_n(\alpha a\beta)\|
=\|T_n[\sum_{i,j=1}^m\alpha_{si}a_{ij}\beta_{jt}]\|
=\|[\sum_{ij=1}^m\alpha_{si}T(a_{ij})\beta_{jt}]\|
=\|\alpha T_m$ $(a)\beta\|=\|\alpha[[D, a_{ij}]]\beta\|
\le\|\alpha\|\|[[D, a_{ij}]]\|\|\beta\|
=\|\alpha\|\|T_m(a)\|\|\beta\|=\|\alpha\|L_m(a)\|\beta\|$. So
$\mathcal{L}$ is a matrix seminorm on $\mathcal{B}$. For any
$[a_{ij}1]\in M_n(\mathbb{C}1)$, we have that $
L_n([a_{ij}1])=\|T_n([a_{ij}1])\|=\|[T(a_{ij}1)]\|=\|[0]\|=0$,
that is, $M_n(\mathbb{C}1)\subseteq K_n$.
\end{proof}

\begin{proposition}\label{pro:42} Let $P_n=M_n(\mathbb{C}1)$ and
$\mathcal{P}=(P_n)$. Then $\mathcal{A}$ has the matrix
2-decomposable property about $\mathcal{P}$.
\end{proposition}

\begin{proof} Given $f\in CB_n^2(\tilde{\mathcal{A}})$. Suppose that $f$ is self-adjoint,
that is, $f^*=f$. By Satz 4.5 in \cite{wi} or Proposition 1.3,
Remark 1.5 and Theorem 1.6 in \cite{ha} and the injectivity of
$M_n$, there are completely positive linear mappings $\psi_1$ and
$\psi_2$ from $\mathcal A$ to $M_n$ such that $f=\psi_1-\psi_2$
and $\|\psi_1+\psi_2\|=\|\psi_1-\psi_2\|_{cb}=\|f\|_{cb}$. But
$0_n=f(1)=\psi_1(1)-\psi_2(1)$. So $\psi_1(1)=\psi_2(1)$.
Therefore,
$\|\psi_1\|=\|\psi_1\|_{cb}=\|\psi_1(1)\|=\|\psi_2(1)\|=\|\psi_2\|_{cb}=\|\psi_2\|$
by Proposition 3.5 in \cite{pau}. Since $\|f\|_{cb}\le 2$, we have
that $\|\psi_1(1)+\psi_2(1)\|=\|\psi_1+\psi_2\|\le 2$. So
$0_n\le\psi_1(1)=\psi_2(1)\le 1_n$. If $\psi_1(1)\le 1_n$ but
$\psi_1(1)\neq 1_n$, choosing a $\phi\in CS_1(\mathcal{A})$ and
setting $\varphi_1(a)=\psi_1(a)+(1_n-\psi_1(1))\phi(a)$ and
$\varphi_2(a)=\psi_2(a)+(1_n-\psi_2(1))\phi(a)$ for $a\in\mathcal
A$, then we have that $\varphi_1(1)=\varphi_2(1)=1_n$ and
$\varphi_1-\varphi_2=\psi_1-\psi_2$. For $k\in\mathbb N$ and
$[a_{ij}]\in M_k(\mathcal{A})$, we have that
$(\varphi_1)_k([a_{ij}])=[\varphi_1(a_{ij})]=[\psi_1(a_{ij})+(1_n-\psi_1(1))\phi(a_{ij})]
=[\psi_1(a_{ij})]+[(1_n-\psi_1(1))\phi(a_{ij})]
=(\psi_1)_k([a_{ij}])+[\phi(a_{ij})]\otimes(1_n-\psi_1(1))
=(\psi_1)_k([a_{ij}])+\phi_k([a_{ij}])\otimes(1_n-\psi_1(1))$.
>From this observation and the completely positivity of $\psi_1$
and $\phi$,  We see that $\varphi_1\in CS_n(\mathcal{A})$.
Similarly, we have that $\varphi_2\in CS_n(\mathcal{A})$. If $f$
is not self-adjoint, we put that $f_1=\frac12(f+f^*),\,
f_2=\frac1{2i}(f-f^*)$. Then $f_1$ and $f_2$ are self-adjoint
completely bounded linear mappings. And if $\|f\|_{cb}\le2$, then
$\|f_1\|_{cb}\le2$ and $\|f_2\|_{cb}\le 2$. So the result above
implies that there exist $\varphi_1, \varphi_2, \varphi_3,
\varphi_4\in CS_n(\mathcal A)$ such that
$f=(\varphi_1-\varphi_2)+i(\varphi_3-\varphi_4)$. Accordingly,
$\mathcal{A}$ has the matrix 2-decomposable property about
$\mathcal{P}$.
\end{proof}

It is apparent that $\mathbb{C}1$ is a closed subspace of
$\mathcal{A}$ and $\mathcal{B}$, so $\mathcal{B}/(\mathbb{C}1)$ is
an operator space with the quotient matrix norm
$\|\cdot\|^\sim=(\|\cdot\|_n^\sim)$ from $\|\cdot\|=(\|\cdot\|_n)$
on $\mathcal{B}$\cite{ru}.

\begin{theorem}\label{th:43} Let $(\mathcal{H}, D)$ be an unbounded Fredholm
module over a unital $C^*$-algebra, and let the matrix metric
$\mathcal{D_L}$ on $\mathcal{CS(A)}$ be defined by
\[\begin{array}{rcl}
D_{L_n}(\varphi, \psi)&=&\sup\big\{\|\ll\varphi, a\gg-\ll\psi,
a\gg\|: a=[a_{ij}]\in M_r(\mathcal{B}),\\
&&L_r(a)=\|[[D, a_{ij}]]\|\le 1, r\in\mathbb{N}\big\}\end{array}\]
for $\varphi, \psi\in CS_n(\mathcal{A})$. Then the
$\mathcal{D_L}$-topology on $\mathcal{CS(A)}$ agrees with the
BW-topology if and only if the image of $L_1^1$ in
$\mathcal{B}/(\mathbb{C}1)$ is totally bounded for
$\|\cdot\|_1^\sim$.
\end{theorem}

\begin{proof} First we show that $K_n=M_n(\mathbb{C}1)$ under the conditions of the theorem.

Suppose that the $\mathcal{D_L}$-topology on $\mathcal{CS(A)}$
agrees with the BW-topology. For each $n\in\mathbb{N}$,
$CS_n(\mathcal{A})$ is BW-compact implies that $D_{L_n}$ is
bounded. If $a=[a_{ij}]\in M_n(\mathcal{B})\setminus
M_n(\mathbb{C}1)$ and $L_n(a)=0$, then there exists an
$a_{ij}\in\mathcal{B}\setminus(\mathbb{C}1)$ such that
$L_1(a_{ij})=0$ since $\mathcal{L}$ is a matrix seminorm by
Proposition \ref{pro:41}. From the definition of $\mathcal{L}$, we
see that $L_n(a^*)=L_n(a)$ for each $a\in M_n(\mathcal{B})$. There
is a self-adjoint $b\in\mathcal{B}\setminus(\mathbb{C}1)$ with
$L_1(b)=0$. So there is a bounded linear functional $f$ on
$\mathcal{A}$ such that $\|f\|=1$, $f(y)=0$ for each $y$ in
$\mathbb{C}1$, and $f(b)\neq 0$. Since $b$ is self-adjoint, there
is a bounded hermitian linear functional $h$ on $\mathcal{A}$ such
that $h(1)=0$ but $h(b)\neq 0$. By Hahn-Jordan decomposition
theorem\cite{ka}, there are positive linear functionals
$\varphi_1$ and $\varphi_2$ on $\mathcal{A}$ such that
$h=\varphi_1-\varphi_2$ and $\|h\|=\|\varphi_1\|+\|\varphi_2\|$. 
From the condition $h(1)=0$, we get $\varphi_1(1)=\varphi_2(1)$.
Thus there are $\psi_1, \psi_2\in CS_1(\mathcal{A})$ and $c>0$
such that $h=c(\psi_1-\psi_2)$. Now, for any $r\in\mathbb{R}^+$ we
have that $L_1(rb)=0<1$. Therefore, $D_{L_1}(\psi_1,
\psi_2)\ge|\psi_1(rb)-\psi_2(rb)|=\frac{r}{c}|h(b)|$ for all
$r\in\mathbb{R}^+$. This contradicts the boundedness of $D_{L_1}$.
So $K_n=M_n(\mathbb{C}1)$ by Proposition \ref{pro:41}.

Assume that the image of $\mathcal{L}^1$ in
$\mathcal{B}/(\mathbb{C}1)$ is totally bounded for
$\|\cdot\|^\sim$. If there is an $a=[a_{ij}]\in K_n$ and $a\notin
M_n(\mathbb{C}1)$, then $\|\tilde{a}\|_n^\sim\neq 0$. Since $ra\in
K_n$ for any $r\in\mathbb{R}$, we get that
$\|\widetilde{ra}\|_n^\sim=|r|\|\tilde{a}\|_n^\sim\to
+\infty(r\to\infty)$. This contradicts the boundedness of the
image of $\mathcal{L}^1$ in $\mathcal{B}/(\mathbb{C}1)$. Thus we
also have that $K_n=M_n(\mathbb{C}1)$ by Proposition \ref{pro:41}.

For $n\in\mathbb N$, let $\varphi=[\varphi_{ij}],
\psi=[\psi_{ij}]\in CS_n(\mathcal{A})$ and $\varphi\neq\psi$. Then
there exist $i_0,j_0\in\{1,2,\cdots,n\}$ such that
$\varphi_{i_0,j_0}\neq\psi_{i_0,j_0}$. Since $\mathcal{B}$ is norm
dense in $\mathcal{A}$ and $\varphi_{i_0,j_0},\psi_{i_0,j_0}\in
\mathcal{A}^*$, there is a $b\in\mathcal{B}$ such that
$\varphi_{i_0,j_0}(b)\neq\psi_{i_0,j_0}(b)$. Therefore,
$\varphi(b)\neq\psi(b)$. So $\mathcal{B}$ separates the points of
$\mathcal{CS(A)}$.

Since a unital and completely contractive linear mapping from
$\mathcal{A}$ into $M_n$ is completely positive(Proposition 3.4 in
\cite{pau}), the conditions in the theorem imply that the matrix
state space $\mathcal{CS(A)}$ of $\mathcal{A}$ is the graded set
$\mathcal{GCS(A)}$ that was defined in section \ref{sec:3} with
$\varphi_n^{(0)}(1)=1_n$. Now the theorem follows from Proposition
\ref{pro:42}, Theorem \ref{th:36} and Theorem \ref{th:38}.
\end{proof}

\begin{remark} Though $D_{L_n}$ does not always define a bounded
metric(see Proposition 4 in \cite{co1}), the conditions in the
theorem above imply that $D_{L_n}$ is bounded from the proof of
it.
\end{remark}

\section{Matrix metrics from matrix Lipschitz seminorms}
\label{sec:5}

A complex vector space $V$ is said to be {\it matrix ordered} if
$V$ is a *-vector space, each $M_n(V)$ is partially ordered and
$\gamma^*M_n(V)^+\gamma\subseteq M_m(V)^+$ when $\gamma\in
M_{n,m}$. A {\it matrix order unit space} $(\mathcal V, 1)$ is a
matrix ordered space $\mathcal V$ together with a distinguished
order unit $1$ satisfying the conditions: $\mathcal V^+$ is a
proper cone with the order unit $1$ and each of the cones
$M_n(\mathcal V)^+$ is Archimedean. Each matrix order unit space
$(\mathcal V, 1)$ may be provided with the norm
$\|v\|=\inf\left\{t\in\mathbb R:\, \left[\begin{array}{cc} t1&v\\
v^*&t1\end{array}\right]\ge 0\right\}$. As in section \ref{sec:3},
we will not assume that $\mathcal V$ is complete for the norm.

A {\it Lipschitz seminorm} on an order unit space $S$ is a
seminorm on $S$ such that its null space is the scalar multiples
of the order unit\cite{rie}. Here we give a non-commutative
version of Lipschitz seminorms.

\begin{definition}\label{def:51} Given a matrix order unit space $(\mathcal{V},
1)$, a {\it matrix Lipschitz seminorm} $\mathcal{L}$ on
$(\mathcal{V},1)$ is a sequence of seminorms $L_n:\,
M_n(\mathcal{V})\longmapsto [0, +\infty)$ such that the null space
of each $L_n$ is $M_n(\mathbb{C}1)$, $L_{m+n}(v\oplus
w)=\max\{L_m(v), L_n(w)\}$, $L_n(\alpha
v\beta)\le\|\alpha\|L_m(v)\|\beta\|$ and $L_m(v^*)=L_m(v)$ for any
$v\in M_m(\mathcal{V})$, $w\in M_n(\mathcal{V})$, $\alpha\in
M_{n,m}$ and $\beta\in M_{m,n}$.
\end{definition}

By an {\it operator system} one means a self-adjoint subspace
$\mathcal{R}$, which is not necessarily closed, in a unital $C^*$-algebra $\mathcal{A}$ such that
$\mathcal{R}$ contains the identity $1\in\mathcal{A}$. Paralled to
the abstract characterization of Kadison's function systems as
order unit space\cite{kad}, Choi-Effros's representation theorem
says that every matrix order unit space is completely order
isomorphic to an operator system\cite{chef}. For a matrix order
unit space $(\mathcal{V}, 1)$, the {\it matrix state space} is the
collection $\mathcal{CS(V)}=(CS_n(\mathcal{V}))$ of matrix states
\[
CS_n(\mathcal{V})=\{\varphi\hbox{: $\varphi$ is a unital
completely positive linear mapping from $\mathcal{V}$ to
$M_n$}\},\] \cite{wewi}. We consider $\mathcal{CS(V)}$ the
non-commutative version of the state space of an order unit space.

\begin{proposition}\label{pro:52} Let $(\mathcal{V}, 1)$ be a matrix order unit
space, and let $P_n=M_n(\mathbb{C}1)$ and $\mathcal{P}=(P_n)$.
Then $\mathcal{V}$ has the matrix 2-decomposable property about
$\mathcal{P}$.
\end{proposition}

\begin{proof}
Since every matrix order unit space is completely order isomorphic
to an operator system, there is a unital $C^*$-algebra
$(\mathcal{A}, 1_{\mathcal{A}})$ such that
$\mathcal{V}\subseteq\mathcal{A}$ and $1_{\mathcal{A}}=1$. Because
$\mathbb{C}1$ is a closed subspace of $\mathcal{V}$, we have the
complete isometry
$(\mathcal{V}/(\mathbb{C}1))^*\cong(\mathbb{C}1)^\perp$, where
$(\mathbb{C}1)^\perp=\{f\in\mathcal{V}: f(1)=0\}$\cite{bl}. So
$(\tilde{\mathcal{V}})^*$ is just the subspace of $\mathcal{V}^*$
consisting of those $f\in\mathcal{V}^*$ such that $f(1)=0$.

For $f\in CB_n^2(\tilde{\mathcal{V}})$, there exists a completely
bounded linear mapping $g$ from $\mathcal{A}$ to $M_n$, which
extends $f$, with $\|g\|_{cb}=\|f\|_{cb}$ by Theorem 7.2 in
\cite{pau}. From Proposition \ref{pro:42}, there are $\psi_1,
\psi_2, \psi_3, \psi_4\in CS_n(\mathcal{A})$ such that
$g=(\psi_1-\psi_2)+i(\psi_3-\psi_4)$. Let
$\varphi_i=\psi_i|_{\mathcal{V}}(i=1,2,3,4)$, the restriction of
$\psi_i$ to $\mathcal{V}$. Clearly each $\varphi_i\in
CS_n(\mathcal{V})$. So we obtain that $
f=(\varphi_1-\varphi_2)+i(\varphi_3-\varphi_4)$. Therefore,
$\mathcal{V}$ has the matrix 2-decomposable property about
$\mathcal{P}$.
\end{proof}

\begin{theorem}\label{th:53}
Let $(\mathcal{V},1)$ be a matrix order unit space with operator
space norm $\|\cdot\|=(\|\cdot\|_n)$, and let $\mathcal{L}=(L_n)$
be a matrix Lipschitz seminorm on $(\mathcal{V},1)$. Then the
$\mathcal{D_L}$-topology on $\mathcal{CS(V)}$ agrees with the
BW-topology if and only if the image of $L_1^1$ in
$\tilde{\mathcal{V}}$ is totally bounded for $\|\cdot\|_1^\sim$.
\end{theorem}

\begin{proof}
If $f$ is a completely bounded linear mapping from $\mathcal{V}$
to $M_n$ and $f(1)=1_n$ and $\|f\|_{cb}=1$, then $f$ is completely
positive by Proposition 3.4 in \cite{pau}. So the matrix state
space $\mathcal{CS(V)}$ of $\mathcal{V}$ is the graded set
$\mathcal{GCS(V)}$ that was defined in section \ref{sec:3} with
$\varphi_n^{(0)}(1)=1_n$. Now by Proposition \ref{pro:52}, Theorem
\ref{th:36} and Theorem \ref{th:38}, the $\mathcal{D_L}$-topology
on $\mathcal{CS(V)}$ agrees with the BW-topology if and only if
the image of $L_1^1$ in $\tilde{\mathcal{V}}$ is totally bounded
for $\|\cdot\|_1^\sim$.
\end{proof}

\end{document}